\documentclass[11pt]{article}
\usepackage{amsfonts}
\usepackage{mathrsfs}
\usepackage{amsfonts,amssymb,mathrsfs}
\usepackage{amsmath,amscd}
\usepackage[dvips]{graphicx}
\usepackage{color}
\usepackage[all]{xy}
\usepackage{graphicx}
\usepackage{mathptm,pslatex}

\definecolor{red}{rgb}{1,0,0}
\definecolor{blue}{rgb}{0,0,1}
\definecolor{darkblue}{rgb}{0,0,0.4}

\begin{document}
\newtheorem{Def}{Definition}[section]
\newtheorem{exe}{Example}[section]
\newtheorem{Prop}[Def]{Proposition}
\newtheorem{Theo}[Def]{Theorem}
\newtheorem{Lem}[Def]{Lemma}
\newtheorem{rem}{\noindent\mbox{Remark}}[section]
\newtheorem{Koro}[Def]{Corollary}
 \newcommand{\Ker}{\rm Ker}
  \newcommand{\Soc}{\rm Soc}
 \newcommand{\lra}{\longrightarrow}
 \newcommand{\ra}{\rightarrow}
 \newcommand{\add}{{\rm add\, }}
\newcommand{\gd}{{\rm gl.dim\, }}
\newcommand{\End}{{\rm End\, }}
\newcommand{\overpr}{$\hfill\square$}
\newcommand{\rad}{{\rm rad\,}}
\newcommand{\soc}{{\rm soc\,}}
\renewcommand{\top}{{\rm top\,}}
\newcommand{\fdim}{{\rm fin.dim}\,}
\newcommand{\fidim}{{\rm fin.inj.dim}\,}
\newcommand{\gldim}{{\rm gl.dim}\,}
\newcommand{\cpx}[1]{#1^{\bullet}}
\newcommand{\D}[1]{{\mathscr D}(#1)}
\newcommand{\Dz}[1]{{\rm D}^+(#1)}
\newcommand{\Df}[1]{{\rm D}^-(#1)}
\newcommand{\Db}[1]{{\mathscr D}^b(#1)}
\newcommand{\C}[1]{{\mathscr C}(#1)}
\newcommand{\Cz}[1]{{\rm C}^+(#1)}
\newcommand{\Cf}[1]{{\rm C}^-(#1)}
\newcommand{\Cb}[1]{{\mathscr C}^b(#1)}
\newcommand{\K}[1]{{\mathscr K}(#1)}
\newcommand{\Kz}[1]{{\rm K}^+(#1)}
\newcommand{\Kf}[1]{{\rm K}^-(#1)}
\newcommand{\Kb}[1]{{\mathscr K}^b(#1)}
\newcommand{\modcat}[1]{#1\mbox{{\rm -mod}}}
\newcommand{\stmodcat}[1]{#1\mbox{{\rm -{\underline{mod}}}}}
\newcommand{\pmodcat}[1]{#1\mbox{{\rm -proj}}}
\newcommand{\imodcat}[1]{#1\mbox{{\rm -inj}}}
\newcommand{\opp}{^{\rm op}}
\newcommand{\otimesL}{\otimes^{\rm\bf L}}
\newcommand{\rHom}{{\rm\bf R}{\rm Hom}\,}
\newcommand{\pd}{{\rm pd}}
\newcommand{\Hom}{{\rm Hom \, }}
\newcommand{\Coker}{{\rm coker}\,\,}
\newcommand{\Ext}{{\rm Ext}}
\newcommand{\im}{{\rm im}}
\newcommand{\Cone}{{\rm Cone}}
\newcommand{\fini}{{finitistic}}
\newcommand{\proj}{{\rm  proj}}
\newcommand{\al}{{\alpha}}
\newcommand{\Tor}{{\rm  Tor}}
\newcommand{\ad}{{\rm  add}}
\newcommand{\StHom}{{\rm \underline{Hom} \, }}
\pagenumbering{arabic}{\Large \bf
\begin{center}
Homological dimensions and strongly idempotent ideals

\end{center}}
\centerline{ {\sc Dengming xu }}
\begin{center}
 {\scriptsize \it  Sino-European Institute of Aviation Engineering,
 Civil Aviation University of China, 300300 Tianjin,
People's Republic of  China\\}

\end{center}
\renewcommand{\thefootnote}{\alph{footnote}}
\setcounter{footnote}{-1} 
\renewcommand{\thefootnote}{\alph{footnote}}
\setcounter{footnote}{-1} 
\renewcommand{\thefootnote}{\alph{footnote}}
\setcounter{footnote}{-1} \footnote{ E-mail:
xudeng17@163.com,\ Tel: 0086-22-24092750}

 \abstract {Let $A$ be an Artin algebra and  $e$  an idempotent in $A$. It is an interesting topic  to compare the homological dimension of the algebras $A,A/AeA$ and $eAe$.
 For example, in \cite{Aus2}, the relation among the global dimension of these algebras  is  discussed  under the condition that $AeA$ is a strongly idempotent ideal. Motivated by this, we try    to   compare the finitistic dimension of these algebras under certain homological conditions on $AeA$. In particular, under the condition that $AeA$ is a strongly idempotent ideal with finite projective dimension,   we prove  that if the finitistic projective (or injective) dimension of  $eAe$ and $A/AeA$ are finite, then the finitistic projective (or injective) dimension of $A$ is finite.
  This is a generalized version of   the main result in \cite{AHLU}.


\section{Introduction}
The open finitistic dimension conjecture says that  the finitistic dimension of any Artin algebra is finite (see \cite{B}). It is known that a positive answer to this conjecture  will imply the solutions to  other  homological conjectures  (see \cite{Y}).
 Thus, it is one of the main topics    in the representation theory of Artin algebras  to study the finiteness of the finitistic dimension of Artin algebras and  it has been proved that  several classes of algebras have finite finitistic dimension. For details,  we refer to  \cite{xu,ZH4} and the references therein.
 Let $A$ be an Artin algebra and $e$ an idempotent in $A$.
In \cite{AHLU},  it is proved that if $e$ is  primitive idempotent in  $A$ such that
 $AeA$ is projective, then the finiteness of the finitistic dimension of $A/AeA$ implies that  of the finitistic dimension of  $A$. Further, the result implies that the finitistic dimension of standardly stratified algebras is  finite.
 On the other hand,
 it is also  an interesting   topic to compare the  homological  dimension of   the three algebras $A,A/AeA$ and $eAe$.
 In general, it is difficult for us  to do it  directly.
 For this reason, the work usually has  been done under  certain homological conditions on the ideal $AeA$  (see \cite{cps,fgr, ful}).  For example, the global dimension of these algebras  are compared under the condition that  $AeA$ is a strongly idempotent ideal in \cite{Aus2}.

Motivated by the statement  above,  in the paper, we try  to  study the relation among the finiteness of the finitistic dimension of $A$, $A/AeA$ and $eAe$ by comparing   the finitistic dimension of these algebras. We shall mainly  do this under the condition that $AeA$ is a strongly idempotent ideal with finite projective dimension. Recall that $AeA$ is called
a strongly idempotent ideal if the epimorphism $A\lra A/AeA$ induces isomorphisms
$\rho_{X,Y}^n: \Ext_{A/AeA}^n(X,Y)\ra \Ext_A^n(X,Y)$ for all $n>0$
 and for all $A/AeA$-modules $X$ and $Y$. A simple example of  strongly idempotent ideals  is the ideal $AeA$ which is projective on either side. Moreover, if we just need that $e$ is an idempotent which is not necessarily  primitive, then we can  construct  examples of strongly idempotent ideal which are not projective but  of  finite projective dimension.  For example,
   let $0=I_0 \subset I_1\subset \cdots \subset I_{n-1}\subset I_n=A  $ be
 a chain of idempotent ideals of $A$ such that  $I_{k+1}/I_k$ is a projective  $A/I_k$-module for $0\leq k\leq n-1$. Then   $I_k$
is a strongly idempotent ideal with finite projective dimension for $1\leq k\leq n$ (see \cite{Aus2}).
 This implies that our consideration is completely new and  more general than that in \cite{AHLU}.
 Denote $\fdim(A)$  the finitistic projective dimension of $A$, $\fidim(A)$  the finitistic injective dimension of $A$,   and $\pd_A(AeA)$ the projective dimension of $AeA$ as a left $A$-module. Now our main result can be  stated as follows.

  \noindent  {\it {\bf Theorem} Let $e$ be an idempotent in  an Artin algebra $A$. Suppose that $AeA$ is a strongly idempotent ideal with finite projective  dimension. Then we have the following.
   \begin{itemize}
     \item [$(1)$]  $\fdim(A)\leq\max\{2\fdim(eAe)+1,\pd_A(AeA)+\fdim (eAe)+\fdim(A/AeA)+2\};$
     \item [$(2)$] $\fidim(A)\leq \pd_{A} (AeA) +\fidim(eAe)+\fidim(A/AeA)+2.$
   \end{itemize}
In particular, if the finitistic projective (or injective) dimension of $eAe$ and $A$ are finite, then the finitistic projective (or injective) dimension   of $A$ is finite. }

The theorem  is a generalized version of   \cite[Theorem 5.4]{Aus2}.
 We will prove it   in Theorems \ref{theo1}  and \ref{theo3}.
Suppose that $e$ is   a primitive idempotent in $A$.  Then $eAe$ is a local algebra. Recall  that the finitistic dimension of  a local Artin algebra is zero. Immediately from the theorem, we get the following result as a  corollary,  which implies that
 both the  finitistic projective and the finitistic injective dimension of standardly stratified algebras are finite (see \cite{AHLU}).

  \noindent  {\it {\bf Corollary} {\rm  \cite{AHLU}} Let $e$ be a primitive idempotent in $A$. Suppose that $AeA$ is  projective. Then we have the following.
   \begin{itemize}
     \item [$(1)$]  $\fdim(A)\leq \fdim(A/AeA)+2.$
     \item [$(2)$] $\fidim(A)\leq \fidim(A/AeA)+2.$
   \end{itemize}
}

The paper is organized as follows. In Section 2, we give some notations, definitions and known results needed in the proof of the theorem. In Section 3, we  prove    basic results in the first subsection; then we give the proof of the main result in the next  two subsections; finally, we prove   general results in the last subsection.
\section{preliminaries}
In this section, we give some notations, definitions and known results  needed in the proof of the theorem.
\subsection{Notations and definitions}
Let $A$ be an Artin algebra. We denote $A^{op}$ the opposite algebra of $A$.  Unless otherwise specified, all modules considered  are finitely generated left $A$-modules. Let $X$ be an $A$-module and $n$ be an integer with $n>0$. The projective dimension of $X$ is denoted by $\pd_A(X)$; the $n$-th syzygy of $X$ is denoted by $\Omega_A^n(X)$. The finitistic projective dimension, or simply the finitistic dimension  of $A$, which is denoted by $\fdim(A)$,  is defined as the supremum  of the projective dimension of finitely generated  $A$-modules  with finite projective dimension. The finitistic injective  dimension of $A$, which is denoted by $\fidim(A)$,  is defined as the supremum  of the injective  dimension of finitely generated  $A$-modules  with finite injective  dimension.  An element $e$ in  $A$ is called an idempotent if $e^2=e$. An  ideal  $I$ in $A$ is called an {\bf idempotent ideal} if  there exists an idempotent $e$  in $A$ such that $I=AeA.$ Now suppose that $e$ is an idempotent in $A$ and write $I=AeA$.
Let $X$ be an $A/I$-module. Then the epimorphism $A\lra A/I$  induces a naturel $A$-module structure of $X$. Let $X$ and $Y$ be two $A/I$-modules. Then  there is an isomorphism $\Hom_{A/I}(X,Y)\simeq\Hom_A(X,Y)$. The isomorphism   induces morphisms
$\rho_{X,Y}^n: \Ext_{A/I}^n(X,Y)\ra \Ext_A^n(X,Y)$  for  all $n>0$ and for all $A/I$-modules $X$ and $Y$.  The ideal $I$ is called a {\bf strongly idempotent ideal} if,  for all $A/I$-modules $X$ and $Y$, $\rho_{X,Y}^n$ are isomorphisms for all $n\geq 0 $ (see \cite{Aus2}).
It is proved in \cite{cps} that  $AeA$ is a strongly idempotent ideal  if and only if
 \begin{itemize}
   \item [(a)] the multiplication map $Ae\otimes_{eAe} eA\lra AeA$ is an isomorphism  and
   \item [(b)]$\Tor_n^{eAe}(Ae,eA)=0$ for  all $n>0$.
 \end{itemize} The ideal satisfying the equivalent condition is  also called a stratifying ideal in \cite{cps}.
  An Artin algebra $A$ is called a {\bf CPS-stratified algebra} if there exists a chain
$$0=I_0 \subset I_1\subset \cdots \subset I_{n-1}\subset I_n=A $$ of idempotent  ideals in $A$
such that $I_{k+1}/I_k$ is a strongly idempotent ideal of $A/I_k$  and  that  $I_{k+1}/I_k$ is generated by a  primitive idempotent in $A/I_k$ for $0\leq k\leq n-1$. In addition, if $I_{k+1}/I_k$ is a projective  $A/I_k$-module for $0\leq k\leq n-1$, then $A$ is called a {\bf standardly stratified algebra}.
 For more  details about the definition, we refer to \cite{cps}.

In the rest of the paper,  we will write $\overline{A}=A/AeA$ and $ B=eAe$ for abbreviation sometimes.

\subsection{Known results}
In this subsection, we collect some results from  \cite{Aus2}, which will be  needed in the proof of our main results.

 Let $e$ be  be an idempotent in $A$. The following lemma gives more characteristics of a strongly idempotent ideal.
\begin{Lem}\label{lem1}{\rm \cite[Propsition 1.3]{Aus2}} The following statements are equivalent.
\begin{itemize}
  \item [$(1)$] The ideal $AeA$ is a strongly idempotent ideal.
  \item [$(2)$] The epimorphism $A\lra A/AeA$ induces isomorphisms
$\Tor^A_n(X,Y)\ra \Tor^{\overline{A}}_n(X,Y)$
for all $\overline{A}$-modules $X$ and $Y$, and  for   all $n>0$.
  \item [$(3)$] $\Tor_n^A(A/AeA, Y)=0$ for all $\overline{A}$-modules $Y$ and all $n>0.$
\end{itemize}
\end{Lem}

 Denote $\ad(Ae)$ the full subcategory of the category of $A$-modules whose objects are  direct summands of   direct sums of  finite copies of $Ae$.  Let $X$ be an $A$-module and
 $$\cdots\lra P_n\lra P_{n-1}\lra \cdots\lra P_0\lra X\lra 0$$
the minimal projective resolution of $X$.
Let $k$ be a naturel number. We say that  $X$ is  in ${\bf P_e^k}$ if  $P_i$ is  in $\ad(Ae)$  for  $0\leq i\leq k$; and say that  $X$ is in ${\bf P_e^\infty}$ if $X$ is in ${\bf P_e^k}$ for all $k\geq 0$. It is known from {\cite[Theorem 2.1]{Aus2}} that  $AeA$ is a strongly idempotent ideal if and only if $AeA$ is in ${\bf P_e^\infty}$. It follows that  if $AeA$ is projective on either side, then it is a strongly idempotent ideal.
 The following lemma provides us  a criteria for determining whether an $A$-module is in  ${\bf P_e^\infty}$.
\begin{Lem}{\label{lem2}}{\rm \cite[Propsition 2.4]{Aus2}} The following conditions are equivalent for an $A$-module $X$.
\begin{itemize}
  \item [$(1)$] The $A$-module $X$ is in ${\bf P_e^\infty}$.
  \item [$(2)$] $\Tor_n^A(A/AeA, X)=0$ for all $n\geq 0.$
\end{itemize}
\end{Lem}

{\it Proof.} For convenience, we include a proof here.  It suffices to show that (2) implies (1).
Let $X$ be an $A$-module such that $\Tor_n^A(A/AeA, X)=0$ for all $n\geq 0.$  Let
$$\cdots\lra P_m\lra P_{m-1}\lra \cdots\lra P_0\lra X\lra 0$$
be the projective resolution of $X$. Since $A/AeA\otimes_AX=0, $ we get that  $P_0$ is  in $\ad(Ae)$.  It follows that $A/AeA\otimes_AP_0=0$. Applying the functor $A/AeA\otimes_A-$ to the exact sequence $0\ra \Omega_A(X)\ra P_0\ra X\ra 0$, we get  an  exact sequence
$$0\lra \Tor_1^A(A/AeA, X)\lra A/AeA\otimes_A\Omega_A(X)\lra A/AeA\otimes_AP_0\lra A/AeA\otimes_A X\lra 0.$$ It follows from $\Tor_1^A(A/AeA, X)=0$ that $A/AeA\otimes_A\Omega_A(X)=0$.
Then we get that  $P_1$ is in $\ad(Ae)$. Now  the result  can be   shown inductively.  $\square$

Finally, we need the following lemma.
\begin{Lem}\label{lem4}{\rm \cite[Corollary 3.2]{Aus2}} Let  $X$   be an $A$-module. If $X$ is  in ${\bf P_e^\infty}$, then  $\pd_A(X)=\pd_B(eX)$.
\end{Lem}

%

\section{Proof of the theorem}

Throughout  this section,  we  denote $A$ an Artin algebra and $e$ an idempotent in $A$, and always assume that $AeA$  is a strongly idempotent ideal.
\subsection{Basic results}
 In this subsection, we prove  basic results which will be used in the  following subsections.

\begin{Lem}\label{lem72} Let $X$ be an $A$-module.  If $\Tor^A_k(A/AeA,X)=0$ for $k\geq 1$, then $AeX$ is in ${\bf P_e^\infty}$.
\end{Lem}

{\it Proof.}    Let $X$ be an $A$-module  such that   $\Tor^A_k(A/AeA,X)=0$ for $k\geq 1$. Since $AeA$ is a strongly idempotent ideal and $X/AeX$ is an $\overline{A}$-module, we get from  Lemma \ref{lem1}
that  $\Tor^A_k(A/AeA, X/AeX)=0$ for $k\geq 1$.
  Applying  the functor $\Hom_A(A/AeA,-)$ to the exact sequence  $$0\ra AeX\ra X\ra X/AeX\ra 0,$$ we get a long exact sequence
  { $$\cdots \lra \Tor^A_{k+1}(A/AeA, X) \lra \Tor^A_{k+1}(A/AeA, X/AeX)\lra \Tor^A_{k}(A/AeA, AeX) \lra \Tor^A_{k}(A/AeA, X)\lra \cdots.$$}
It follows that $ \Tor^A_{k}(A/AeA, AeX)=0$ for $k\geq 1.$   Note that  $(A/AeA)\otimes_A AeX=0$. Consequently,
   $ \Tor^A_{k}(A/AeA, AeX)=0$ for $k\geq 0.$ It follows from Lemma \ref{lem2}
   that $AeX$ is  in ${\bf P_e^\infty}$. $\square$

Now we can prove the following result, which will be frequently used in the proof of the main result.
\begin{Lem}\label{lem71} Let $X$ be an $A$-module. Suppose  that there exists  a  naturel number $n$ such that $\Tor_k^B(Ae,eX)=0$ for $k\geq n+1.$  Then $Ae\Omega_A^{n+1}(X)$
is in ${\bf P_e^\infty}$.
\end{Lem}

{\it Proof.}  Denote $D$ the usual duality for Artin algebras. Let $Z$ be an $A^{op}$-module and $Y$ an $A$-module. It is known that there is an isomorphism  $\Tor_k^A(Z,Y)\simeq D\Ext_{A^{op}}^k(Z,D(Y))$ for  each $k\geq 0$.  Let $X$ be an $A$-module  and $m$ a naturel number.
Then we have isomorphisms
$$\begin{array}{ll}\Tor_m^B(Ae,eX)&\simeq D\Ext_{B^{op}}^m(Ae, D(eX))\\ &\simeq D\Ext_{B^{op}}^m (\Hom_{A^{op}}(eA,AeA), \Hom_{A^{op}}(eA,D(X)))\\ &\simeq D\Ext_{A^{op}}^m(AeA,D(X))\\&\simeq \Tor_m^A(AeA,X),
\end{array}$$
where the third isomorphism follows from \cite[Theorem 3.2]{Aus2} and the fact that $AeA$ is in ${\bf P_e^\infty}$.  By assumption,  there exists  a  naturel number $n$ such that $\Tor_k^B(Ae,eX)=0$ for $k\geq n+1.$
Then   $ \Tor_k^A(AeA,X)=0$ for  $k\geq n+1.$
It follows that $ \Tor_k^A(A/AeA,X)=0$ for  $k\geq n+2.$ Therefore,  $ \Tor_k^A(A/AeA,\Omega_A^{n+1}(X))=0$ for  $k\geq 1.$ It follows from Lemma \ref{lem72} that  $Ae\Omega_A^{n+1}(X)$
is in ${\bf P_e^\infty}$. $\square$

 \subsection{Finitistic projective dimensions}
In this subsection, we always assume that $AeA$ is a strongly idempotent ideal with finite projective dimension. This subsection aims to  prove that the finiteness of
   the finitistic dimension of $A/AeA$ and $eAe$ implies  that  of the finitistic dimension of $A$ under this condition.

We need the following homological fact first. For completeness, we include a proof here.
\begin{Lem}\label{lem5}  Let $X $ be an $A$-module. If $\pd_A(X)<+\infty,$ then $\pd_B(eX)< +\infty$
\end{Lem}

 {\it Proof.}  Let $X$ be an $A$-module with finite projective dimension. Let
$$0\lra P_m\lra P_{m-1}\lra \cdots\lra P_0\lra X\lra 0$$
be the projective resolution of $X$. Applying the functor $\Hom(Ae,-)$ to it, we get a long exact sequence
$$0\lra eP_m\lra eP_{m-1}\lra \cdots\lra eP_0\lra eX\lra 0$$ of $B$-modules. Since the strongly idempotent ideal $AeA$   is in ${\bf P_e^\infty}$, we know from Lemma \ref{lem4} that $\pd_B(eA)=\pd_A(AeA)$. Since $\pd_A(AeA)<+\infty$, we have  $\pd_B(eA)<+\infty$.  It follows from the previous long exact sequence that $\pd_B(eX)<+\infty.$ $\square$

\begin{rem}{\rm \label{rem1} If we drop the condition that $AeA$ has finite projective dimension, then the lemma does not have to be true. For example,
let  $A$ be the algebra given by the following quiver
with relation.
\begin{center}
\setlength{\unitlength}{1mm}
\begin{picture}(40,2)(0,0)
\thinlines \put(-1,1){{1}}
 \put(16,1){{2}}
\put(1,0){\circle*{1}} \put(3,-1){\vector(1,0){10}}
\put(13,1){\vector(-1,0){10}} \put(15,0){\circle*{1}}
\put(8,2){$\scriptscriptstyle \beta$} \put(8,-4){$\scriptscriptstyle
\alpha$} \put(28,0){$\scriptstyle {\alpha\beta\alpha=0 }$.}
\end{picture}
\end{center}

Then the indecomposable projective modules of $A$ are as follows:
$$ \begin{array}{c}
1\\ 2\\1\\2\end{array}\qquad \quad \begin{array}{c}
 2\\1\\2\end{array} $$
Let  $S$ be the simple $A$-module corresponding  to   the vertex $1$ and $e_1$ the idempotent corresponding  to   the vertex $1$. It is easy to check that $Ae_1A$ is in ${\bf P_e^\infty}$. It follows that $Ae_1A$ is a strongly idempotent ideal with infinite projective dimension.  One can check that $\pd_A(S)=1$ but $\pd_B(e_1S)$ is infinite.
}
\end{rem}

\begin{Lem}\label{lem7} Let $X$ be an $A$-module with finite projective dimension.   Then there exists a natural number $n$ with  $n\leq \fdim (B)$ such that $Ae\Omega_A^{n+1}(X)$ is in ${\bf P_e^\infty}$.
\end{Lem}

{\it Proof.}    Let $X$ be an $A$-module with finite projective dimension. Then we get  from Lemma \ref{lem5} that  the projective dimension of $eX$ is finite. Suppose that  $\pd_B(eX)=n$. Then  $\Tor^B_k(Ae, eX)=0$ for $k\geq n+1$. Since $AeA$ is a strongly idempotent ideal, we get from Lemma \ref{lem71}
   that $Ae\Omega_A^{n+1}(X)$ is  in ${\bf P_e^\infty}$. $\square$

Now we can prove the first main result in the paper. Although we can also prove the next theorem by Lemma  \ref{lem9} in the next subsection, we include a different proof here, not only since the upper bound here is better but also because  the proof may be  of its own interest.

\begin{Theo}\label{theo1}  Let $e$ be an idempotent  in $A$ such that $AeA$ is a strongly idempotent ideal with finite projective dimension. Then we have the following.
 \begin{itemize}
   \item [$(1)$] $ \fdim(A)\leq\max\{2\fdim(eAe)+1,\pd_A(AeA)+\fdim (eAe)+\fdim(A/AeA)+2\}$
   \item [$(2)$] $\fdim(A)\leq 2\fdim(eAe)+\fdim(A/AeA)+2.$
 \end{itemize}
 In particular,  if  $\fdim(A/AeA)<+\infty$ and $\fdim(eAe)<+\infty$, then   $\fdim(A)<+\infty$.
\end{Theo}

{\it Proof.}
 (1) \ Let  $X$ be an $A$-module with finite projective dimension. We know from Lemma \ref{lem7}  that there exists a natural number $n$ with  $n\leq \fdim (B)$ such that $Ae\Omega_A^{n+1}(X)$ is in ${\bf P_e^\infty}$.  It follows from Lemma \ref{lem4} that $\pd_A(Ae\Omega_A^{n+1}(X))=\pd_B(e\Omega_A^{n+1}(X))$. Since $\pd_A(\Omega_A^{n+1}(X))<+\infty$,  we get from Lemma \ref{lem5}
 that $\pd_B (e\Omega_A^{n+1}(X))< +\infty$.  It follows that $\pd_A(Ae\Omega_A^{n+1}(X))< +\infty$. As a result, $\pd_A(\Omega_A^{n+1}(X)/Ae\Omega_A^{n+1}(X))<+\infty.$ Since $AeA$ is a strongly idempotent ideal, we know that $\pd_{\overline{A}}(Y)\leq \pd_A(Y)$ for any   $\overline{A}$-module $Y$.
  It follows that $\pd_{\overline{A}}(\Omega_A^{n+1}(X)/Ae\Omega_A^{n+1}(X))$ is finite. On the other hand, we get from the change of rings theorem  that $\pd_A(Y)\leq \pd_{\overline{A}}(Y)+\pd_A(\overline{A})$ for any   $\overline{A}$-module $Y$. Then we have
$$\begin{array}{ll}\pd_{{A}}(\Omega_A^{n+1}(X)/Ae\Omega_A^{n+1}(X))&\leq \pd_{\overline{A}}(\Omega_A^{n+1}(X)/Ae\Omega_A^{n+1}(X))+\pd_A(A/AeA)\\
  &\leq   \pd_{\overline{A}}(\Omega_A^{n+1}(X)/Ae\Omega_A^{n+1}(X))+\pd_A(AeA)+1
  \end{array}$$
   \noindent
 Considering the exact sequence $0\ra Ae\Omega_A^{n+1}(X)\ra \Omega_A^{n+1}(X)\ra\Omega_A^{n+1}(X)/Ae\Omega_A^{n+1}(X)\ra 0, $ we get that
 $$\begin{array}{ll}
 \pd_A(X)&\leq n+1+\pd_A(\Omega_A^{n+1}(X))\\
 &\leq n+1+\max\{\pd_A( Ae\Omega_A^{n+1}(X)), \pd_A(\Omega_A^{n+1}(X)/Ae\Omega_A^{n+1}(X))\}\\
 &\leq n+1+\max\{\pd_A( Ae\Omega_A^{n+1}(X)), \pd_{\overline{A}}(\Omega_A^{n+1}(X)/Ae\Omega_A^{n+1}(X))+\pd_A(AeA)+1\}\\
&=  n+1+\max\{\pd_B( e\Omega_A^{n+1}(X)), \pd_{\overline{A}}(\Omega_A^{n+1}(X)/Ae\Omega_A^{n+1}(X))+\pd_A(AeA)+1\}\\
&\leq  \fdim (eAe)+1+\max\{\fdim(eAe), \fdim(A/AeA)+\pd_A(AeA)+1\}\\
&=\max\{2\fdim(eAe)+1,\pd_A(AeA)+\fdim (eAe)+\fdim(A/AeA)+2\}
  \end{array}$$
where the equality follows from  Lemma \ref{lem4} and  the fact that $Ae\Omega_A^{n+1}(X)$  is  in ${\bf P_e^\infty}$. Consequently, $$\fdim(A)\leq\max\{2\fdim(eAe)+1,\pd_A(AeA)+\fdim (eAe)+\fdim(A/AeA)+2\}.$$

(2) Since  the strongly idempotent ideal $AeA$ is in ${\bf P_e^\infty}$. We get  from { Lemma \ref{lem4} that $\pd_A(AeA)=\pd_{eAe}(eA)$.  Then the result follows from  (1).
$\square$

%
\noindent \begin{rem} {\rm
(1) Let  $0=I_0 \subset I_1\subset \cdots \subset I_{n-1}\subset I_n=A $ be a chain of idempotent  ideals in $A$
such that $I_{k+1}/I_k$ is a projective  $A/I_k$-module for $0\leq k\leq n-1$.  It  is proved in  \cite[Proposition 6.1]{Aus2} that  $I_k$ is a strongly idempotent ideal with finite projective dimension  for $1\leq k\leq n.$  This can  provide us   examples of strongly idempotent ideals which are not projective but of  finite projective dimension.

(2) Let $d$ be the supremum of the set $\{\pd_A(X)\mid eX=0\}$. In \cite{ful}, it is proved that $\fdim(A)\leq \fdim(eAe)+d+1$ under the condition that $\pd_{eAe}eA<+\infty$.   This is different from our result since   the supremum  does not have to be  finite under the condition in the theorem.

}\end{rem}

Recall  that if $AeA$ is projective, then  it is a strongly idempotent ideal. Combining with Theorem \ref{theo1}, we get the following result, which implies that the finitistic dimension of standardly stratified algebras is finite.
\begin{Koro}{\rm \cite[Theorem 2.2]{AHLU}}\label{coro1}  Let $e$ be a primitive idempotent in $A$. Suppose that $AeA$ is projective. Then $\fdim(A)\leq \fdim(A/AeA)+2.$
\end{Koro}

 {\it Proof.} Since $e$ is a primitive idempotent in $A$, we know that $eAe$ is a local algebra. Then the result  follows from Theorem \ref{theo1} and the fact that the finitistic dimension of a local Artin algebra is zero. $\square$

\subsection{Finitistic injective dimensions}
This subsection is devoted to showing that the   finiteness of
   the finitistic  injective dimension of $A/AeA$ and $eAe$ implies that  of the finitistic injective  dimension of $A$ under the condition that $AeA$ has    finite projective dimension.

We prove the following lemma first.

\begin{Lem}{\label{lem9}}
Let $0\lra X\lra Y\lra Z\lra 0 $ be an exact sequence  of $A$-modules with $Z$ an $A/AeA$-module. Suppose that $X$ is in ${\bf P_e^\infty}$ and that the projective dimension of $Y$ is finite. Then we have the following.
\begin{itemize}
  \item [$(1)$] $\pd_{\overline{A}}(Z)<+\infty$.
  \item [$(2)$] $\pd_A(Y)\leq \fdim(A/AeA)+\fdim(eAe)+1.$
\end{itemize}
\end{Lem}

{\it Proof.} (1) Let $0\lra X\lra Y\lra Z\lra 0 $ be an exact sequence  of $A$-modules. Suppose that $\pd_A(Y)=n<+\infty.$  Then we get an exact sequence $0\lra \Omega_A^n(X)\lra \Omega_A^n(Y)\oplus Q\lra \Omega_A^n(Z)\lra 0$ with  $\Omega_A^n(Y)\oplus Q$ an  projective $A$-module. It follows that $\Omega^{n+2}_A(Z)\simeq \Omega_A^{n+1}(X)$.  Since $X$ is in  ${\bf P_e^\infty}$, we get that $\Omega_A^{n+1}(X)$ is in ${\bf P_e^\infty}$. It follows that $\Omega^{n+2}_A(Z)$ is in  ${\bf P_e^\infty}$. Suppose that $Z$ is an $A/AeA$-module.
 Let $$\cdots\lra P_{n+2}\lra P_{n+1}\lra \cdots \lra P_0\lra Z\lra 0$$ be the minimal projective resolution of $Z$ as an $A$-module. Since  $AeA$ is a strongly idempotent ideal, we get  from Lemma \ref{lem1} that
 $\Tor_k^A(A/AeA, Z)=0$ for $k\geq 0.$ Applying the functor $A/AeA\otimes _A-$ to the resolution,
we get a long exact sequence
$$\cdots\lra A/AeA\otimes _AP_{n+2}\lra A/AeA\otimes _A P_{n+1}\lra \cdots\lra A/AeA\otimes _A P_0\lra Z\lra 0,$$  which is  the minimal projective resolution of $Z$ as an $A/AeA$-module  since the functor $A/AeA\otimes_A-$ preserves projective covers. We have seen that  $\Omega^{n+2}_A(Z)$ is in  ${\bf P_e^\infty}$. It follows that  $A/AeA\otimes _AP_{n+2}=0$. Therefore, $\pd_{\overline{A}}(Z)\leq n+1<+\infty.$

(2) It follows from (1) that $\pd_{\overline{A}}(Z)<+\infty.$ Suppose $\pd_{\overline{A}}(Z)=m.$
 Comparing the long exact sequences in the proof of (1), we know that $A/AeA\otimes_A P_k=0$ for $k\geq m+1.$ It follows that $\Omega_{A}^{m+1}(Z)$ is in ${\bf P_e^\infty}.$ By assumption, $X$  is in ${\bf P_e^\infty}.$  Considering  the exact sequence $0\ra \Omega_{A}^{m+1}(X)\ra \Omega_{A}^{m+1}(Y)\oplus R\ra \Omega_{A}^{m+1}(Z)\ra 0$ with $R$ a projective $A$-module, we get  that $\Omega_{A}^{m+1}(Y)$ is in ${\bf P_e^\infty}$ since  ${\bf P_e^\infty}$ is closed under extensions. Then we have by Lemma \ref{lem4} that $\pd_A(\Omega_{A}^{m+1}(Y))=\pd_B(e\Omega_{A}^{m+1}(Y))$. It follows that
 $$\pd_A(Y)\leq m+1+\pd_A(\Omega_{A}^{m+1}(Y))=m+1+\pd_B(e\Omega_{A}^{m+1}(Y))\leq \fdim(A/AeA)+\fdim(eAe)+1.$$
 This finishes the proof of the lemma. $\square$

Immediately from the lemma, we get  the following result.

\begin{Prop}  Suppose that   $AeA$ is   a strongly idempotent ideal. Denote $\fdim_A(A/AeA)$ the supremum of the set
$\{\pd_A(X)\mid eX=0 \ \mbox{and}\  \pd_A(X)< +\infty\}$. Then $\fdim_A(A/AeA)\leq \fdim(A/AeA)+\fdim(eAe)+1.$
\end{Prop}

{\it Proof.} Let $X$ be an $\overline{A}$-module. Then $AeX=0$.  It follows from  Lemma \ref{lem9} that
$$\pd_A(X)\leq  \fdim(A/AeA)+\fdim(eAe)+1.$$ This finishes the proof. $\square$


\begin{Lem}\label{lem8} Suppose that $AeA$ is a strongly idempotent ideal with   finite projective dimension as a right $A$-module.  Suppose that $\pd_{A^{op}}(AeA)=n.$ Let $X$ be an $A$-module. Then  $Ae \Omega_A^{n+1}(X)$ is in ${\bf P_e^\infty}$.
\end{Lem}

{\it Proof.} Let $X$ be an $A$-module. Suppose that  $\pd_{A^{op}} (AeA)=n$,  we get from Lemma \ref{lem4}  that $\pd_{B^{op}}(Ae)=n$.  Thus $\Tor_k^B(Ae, eX)=0$ for $k\geq n+1.$
 It follows  from Lemma \ref{lem71}
   that $Ae\Omega_A^{n+1}(X)$ is  in ${\bf P_e^\infty}$. $\square$

Now we can prove the following result, which  is a generalized version of  \cite[Theorem5.4]{Aus2}.
 \begin{Prop} \label{theo2}  Suppose that $AeA$ is a strongly idempotent ideal with  finite projective  dimension as a right $A$-module. Then $$\fdim(A)\leq \pd_{A^{op}}(AeA) +\fdim(eAe)+ \fdim(A/AeA)+2.$$  In particular,  if  $\fdim(A/AeA)<+\infty$ and $\fdim(eAe)<+\infty$, then   $\fdim(A)<+\infty$.
 \end{Prop}

 {\it Proof.}    Let $X$ be an $A$-module with finite projective dimension.  Suppose that  $n=\pd_{A^{op}}(AeA)$. It follows from Lemma \ref{lem8} that  $Ae\Omega_A^{n+1}(X)$ is  in ${\bf P_e^\infty}$.  Then the exact sequence $$0\lra Ae\Omega_A^{n+1}(X)\lra \Omega_A^{n+1}(X)\lra \Omega_A^{n+1}(X)/Ae\Omega_A^{n+1}(X)\lra 0$$ satisfies all the conditions in Lemma \ref{lem9}.  It follows  that $\pd_A(\Omega_A^{n+1}(X))\leq   \fdim(A/AeA)+\fdim(eAe)+1.$ Then we have
 $$\pd_A(X)\leq n+1+\pd_A(\Omega_A^{n+1}(X))\leq n +\fdim(eAe)+ \fdim(A/AeA)+2.$$
 Consequently, $\fdim(A)\leq  \pd_{A^{op}}(AeA)+\fdim(eAe)+ \fdim(A/AeA)+2.\quad \square$

\begin{rem} {\rm Let $A$ be the algebra given    in Remark \ref{rem1}. Let $e_2$ be the idempotent corresponding to the vertex $2$.   Then we  see that  $Ae_2A$ is projective as a left $A$-module. But the projective dimension of $Ae_2A$ as a right $A$-module is infinite.  This shows  that the  previous  proposition can not be deduced from Theorem \ref{theo1}.}
\end{rem}

Now we can prove the following.
 \begin{Theo} \label{theo3} Suppose that $AeA$ is a strongly idempotent ideal with finite projective dimension. Then
 $$\fidim(A)\leq \pd_{A} (AeA) +\fidim(eAe)+\fidim(A/AeA)+2.$$  In particular,  if  $\fidim(A/AeA)<+\infty$ and $\fidim(eAe)<+\infty$, then   $\fidim(A)<+\infty$.
 \end{Theo}

{\it Proof.}  Suppose that $AeA$ is a strongly idempotent ideal with finite projective dimension. Then $\pd_{(A^{op})^{op}}(AeA)$ is finite.   It follows from Proposition \ref{theo2}  that
$$\fdim(A^{op})\leq \pd_{A} (AeA) +\fdim((eAe)^{op})+\fdim((A/AeA)^{op})+2.$$
 Then the result follows from the fact that $\fdim(\Lambda^{op})=\fidim(\Lambda)$ for any Artin algebra $\Lambda$.
 $\square$

 Immediately from the previous theorem, similarly as in Corollary \ref{coro1}, we get the following result, which implies that the finitistic injective dimension of standardly stratified algebras is finite (see \cite{AHLU}).

 \begin{Koro} Let $e$ be a primitive idempotent in $A$. Suppose that $AeA$ is  projective. Then $$\fidim(A)\leq \fidim(A/AeA)+2.$$
 \end{Koro}

 \subsection{General cases }
 In this subsection,  we drop  the condition that $AeA$ has finite projective dimension and strengthen other conditions.
 We  prove  the following  general result first.
 \begin{Prop} \label{prop2} Suppose   that $AeA$ is a strongly idempotent ideal and that  there exists a naturel number $n$ such that $\Tor^B_k(Ae,eX)=0$ for all $A$-modules $X$ with finite projective dimension and all $k\geq n+1$.  Then $$\fdim(A)\leq n +\fdim(eAe)+\fdim(A/AeA)+2.$$  In particular, if  $\fdim(A/AeA)<+\infty$ and $\fdim(eAe)<+\infty$, then   $\fdim(A)<+\infty$.
 \end{Prop}

 {\it Proof.} Let $X$ be an $A$-module with finite projective dimension. By assumption,  $\Tor^B_k(Ae,eX)=0$ for   $k\geq n+1$.  It follows from Lemma \ref{lem71} that $Ae\Omega_A^{n+1}(X)$ is in $ {\bf P_e^\infty}.$  Considering the exact sequence$$0\lra Ae\Omega_A^{n+1}(X)\lra \Omega_A^{n+1}(X)\lra \Omega_A^{n+1}(X)/Ae\Omega_A^{n+1}(X)\lra 0,$$ we get from Lemma \ref{lem9} that
 $$\pd_A(X)\leq n+1+\pd_A(\Omega_A^{n+1}(X))\leq   n +\fdim(eAe)+\fdim(A/AeA)+2.$$
 Thus, $\fdim(A)\leq n +\fdim(eAe)+\fdim(A/AeA)+2.$  Consequently, if  $\fdim(A/AeA)<+\infty$ and \\ $\fdim(eAe)<+\infty$, then   $\fdim(A)<+\infty$. $\square$

Let $A$ be an Artin algebra and $n$ an integer with $n\geq 1$. Set
\begin{center}$\Omega_A^n(A$-\!\!\!\!$\mod):=\{\Omega_A^n(X)\mid X\in
A$-\!\!\!\!$\mod\}$.\end{center} Then  we say that $A$ is   {\bf $*$-syzygy
finite} if there exists a naturel number $m$ such that there exist only  finitely many non-isomorphic indecomposable
 modules  in  $\Omega_A^m(A$-\!\!\!\!$\mod)$.
Recall that an Artin algebra $A$ is said to be of {\bf finite representation type} if there exist only finitely many non-isomorphic indecomposable $A$-modules. Thus, an algebra which is  representation type is $*$-syzygy finite. Obviously,  if $A$
 is $*$-syzygy finite, then $\fdim(A)<+\infty.$

\begin{Prop} \label{prop1}Suppose that $AeA$ is a strongly idempotent ideal such that  $eAe$ is $*$-syzygy finite. If $\fdim(A/AeA)<+\infty,$ then $\fdim(A)<+\infty.$

\end{Prop}

{\it Proof.} Suppose that $eAe$ is
$*$-syzygy finite. Then there exists a naturel number $n$ such that there are only  finitely many non-isomorphic indecomposable
 modules  in  $\Omega_B^n(B$-\!\!\!\!$\mod)$.
  Set
\begin{center}$\mathcal{M}=\{Y\in\Omega_B^n(B$-\!\!\!\!$\mod)\mid $ there  exists  a naturel number $l$ such that  $\Tor_k^B(Ae,Y)=0$ for $k\geq l+1$ $\}.$
\end{center}
 Then there exists a naturel number $p$ such that $\Tor_k^B(Ae,Y)=0$ for all $Y\in \mathcal{M}$ and for all $k\geq p+1.$
Let $X$ be an $A$-module  with $\pd_A(X)=m\geq n+2$. Let
$$0\lra P_m\lra P_{m-1}\lra \cdots\lra P_0\lra X\lra 0$$
be  the  minimal projective resolution of $X$. Applying the functor $\Hom(Ae,-)$ to it, we get a long exact sequence
$$0\lra eP_m\lra eP_{m-1}\lra \cdots\lra eP_0\lra eX\lra 0$$ of $B$-modules.
  Since $AeA$ is a strongly idempotent ideal, we get that $\Tor_k^B(Ae,eA)=0$ for $k\geq 1.$ Thus,
   $\Tor_k^B(Ae,eP_m)=\Tor_k^B(Ae,eP_{m-1})=0$ for $k\geq 1.$
It follows that $\Tor_k^B(Ae,e\Omega_{A}^{m-1}(X))=0$ for $k\geq 2.$ Inductively, one can prove that $\Tor_k^B(Ae,eX)=0$ for $k\geq m+1.$  Then we have $\Tor_k^B(Ae,\Omega_B^{n}(eX))=0$ for $k\geq m-n+1.$  It follows that $\Omega_B^{n}(eX)\in\mathcal{M}$.
Consequently, $\Tor_k^B(Ae,\Omega_B^{n}(eX) )=0$ for $k\geq p+1.$ As a result, $\Tor_k^B(Ae,eX )=0$ for $k\geq n+p+1.$
Note that $\fdim(eAe)<+\infty$ since it is $*$-syzygy finite. It follows from  Proposition \ref{prop2} that  if $\fdim(A/AeA)<+\infty,$ then $\fdim(A)<+\infty.$
 $\square$

 As a direct consequence of the proposition, we get the following.
 \begin{Koro} Suppose that   $AeA$ is   a strongly idempotent ideal and  that   the finitistic dimension of $A/AeA$ is finite. Then the finitistic dimension of $A$ is finite if one of the following conditions holds.
 \begin{itemize}
   \item [$(1)$] The global dimension of $eAe$ is finite.
   \item [$(2)$] The algebra $eAe$ is of finite representation type.
   \item [$(3)$] The algebra $eAe$ is   monomial.
    \end{itemize}

 \end{Koro}

{\it Proof.}  By Proposition \ref{prop1},  it suffices to show that $eAe$ is $*$-syzygy finite under each condition. This is clear for (1) and (2) and it follows from \cite[Theorem I]{hz} that a monomial algebra is  $*$-syzygy finite. This completes the proof. $\square$

Recall that an $A$-module $X$ is said to be {\bf  Gorenstein projective } if there exists an exact sequence $${ P^\bullet}:\cdots \lra P_1\lra P_0\lra  P^0\ra P^1\lra\cdots$$ of projective modules with $X\simeq {\rm Im}(P_0\lra  P^0)$ such that $\Hom_A({P^\bullet}, A)$ is exact for all projective modules $Q$. Let $n$ be a naturel number.  We say that the Gorenstein projective dimension of $X$ is at most $n$ if  $X$ has a Gorenstein projective resolution of length $n$. For more details on  the definition, we refer to \cite{h}.

\begin{Prop} Suppose that $AeA$ is a strongly idempotent ideal such that the Gorenstein projective dimension of  $AeA_A$ is finite. If  $\fdim(A/AeA)<+\infty$ and $\fdim(eAe)<+\infty$, then   $\fdim(A)<+\infty$.
\end{Prop}

{\it Proof.}  Suppose that the Gorenstein  projective  dimension of $AeA_A$ is at most  $n$. Let $X$ be an $A$-module with finite projective dimension.
Then by \cite[Lemma 4.1]{xi2}, we know that $\Tor ^A_k(AeA,X)=0$ for $k\geq n+1$. It follows that
$\Tor ^A_k(A/AeA,X)=0$ for $k\geq n+2$, and therefore $\Tor ^A_k(A/AeA,\Omega^{n+1}_A(X))=0$ for $k\geq 1$. It follows from Lemma \ref{lem72} that $Ae\Omega_{A}^{n+1}(X)$ is in ${\bf P_e^\infty}$. Considering the exact sequence$$0\lra Ae\Omega_A^{n+1}(X)\lra \Omega_A^{n+1}(X)\lra \Omega_A^{n+1}(X)/Ae\Omega_A^{n+1}(X)\lra 0,$$ we get from Lemma \ref{lem9} that
 $$\pd_A(X)\leq n+1+\pd_A(\Omega_A^{n+1}(X)) \leq   n +\fdim(eAe)+\fdim(A/AeA)+2.$$
 Then the result follows from that assumption that  the finitistic  dimension of $eAe$ and $A/AeA$ are finite. $\square$

\medskip

%
%
%
%
\bigskip

{\bf Acknowledgement} \  The work is  supported by the  scientific research
foundation of Civil Aviation University of China (No.2010QD09X).


2013-1-4

\end{document}